\documentclass[10pt]{article}
\usepackage{amssymb}
\usepackage{amsmath}
\usepackage[dvips]{graphics}
\usepackage{epsfig}

\begin{document}
\title{Examples of fourth-order scattering-type operators with embedded eigenvalues in their continuous spectra}

\author{
\textbf{Vassilis G. Papanicolaou}
\\\\
Department of Mathematics
\\
National Technical University of Athens
\\
Zografou Campus, 157 80, Athens, GREECE
\\
\underline{\tt papanico@math.ntua.gr}}
\maketitle

\begin{abstract}
We give examples of fourth-order scattering-type operators, acting on $L_2(\mathbb{R})$, which have eigenvalues embedded in their continuous spectra.
\end{abstract}

\\


\medskip

Before we start, let us mention that an excellent reference for the general spectral theory of linear differential operators, which the reader may want to consult, is Naimark's monograph \cite{N}.

In the first two examples the fourth-order operator $L$ has the form $L = \frac{d^4}{dx^4} + q(x)$ where the ``potential" $q(x)$ has compact support.

\section{An example with a singular potential}
Consider the function
\begin{equation}
\theta(x) =
\left\{
  \begin{array}{ccc}
    \sin x, & \quad &\ 0 \leq x < 3\pi/4; \\
    \\
    A e^{-x}, & \quad &\ 3\pi/4 \leq x.  \\
  \end{array}
\right.
\label{A1}
\end{equation}
We choose
\begin{equation}
A = \frac{1}{\sqrt{2}}\, e^{3\pi/4},
\label{A2}
\end{equation}
so that $\theta(x) \in C^1[0, \infty)$. In particular,
\begin{equation}
\theta\left(\frac{3\pi}{4}\right) = \frac{1}{\sqrt{2}}
\qquad \text{and} \qquad
\theta'\left(\frac{3\pi}{4}\right) = -\frac{1}{\sqrt{2}},
\label{A3}
\end{equation}
while $\theta''(x)$ has a jump discontinuity at $x = 3\pi/4$.

Differentiating $\theta(x)$ four times we get
\begin{equation}
\theta''''(x) - \sqrt{2} \,\, \delta'\left(x - \frac{3\pi}{4}\right) + \sqrt{2} \,\, \delta\left(x - \frac{3\pi}{4}\right) = \theta(x),
\label{A4}
\end{equation}
where $\delta(\cdot)$ is the Dirac delta function.

\medskip

From the basic distribution theory we have that:

(i) $(x - b)\, \delta'(x - b) = - \delta(x - b)$

and

(ii) if $f, g \in C^1(\mathbb{R})$, then $f(b) = g(b)$ implies the distributional equality
\begin{equation*}
f(x) \, \delta(x - b) = g(x) \, \delta(x - b),
\end{equation*}
while the two equalities $f(b) = g(b)$ and $f'(b) = g'(b)$ imply the distributional equality
\begin{equation*}
f(x) \, \delta'(x - b) = g(x) \, \delta'(x - b).
\end{equation*}

\medskip

Using the above properties in \eqref{A4} we obtain
\begin{equation}
\theta''''(x) - 2 \,\, \delta'\left(x - \frac{3\pi}{4}\right) \theta(x) + 4 \,\, \delta\left(x - \frac{3\pi}{4}\right) \theta(x) = \theta(x).
\label{A5}
\end{equation}

Now, let $f(x)$ be the odd extension of $\theta(x)$ on $\mathbb{R}$. Clearly $f(x) \in C^1(\mathbb{R}) \cap L_2(\mathbb{R})$ and \eqref{A5} yields
\begin{equation}
f''''(x) + q_s(x) f(x) = f(x),
\qquad
x \in \mathbb{R},
\label{A6}
\end{equation}
where
\begin{equation}
q_s(x) = Q_s(x) - Q_s(-x),
\label{A7a}
\end{equation}
with
\begin{equation}
Q_s(x) = - 2 \,\, \delta'\left(x - \frac{3\pi}{4}\right) + 4 \,\, \delta\left(x - \frac{3\pi}{4}\right),
\qquad
x \in \mathbb{R}.
\label{A7}
\end{equation}

Therefore, $f(x)$ is an $L_2(\mathbb{R})$-eigenfunction of the scattering-type operator
\begin{equation}
Lu = u'''' + q_s(x) u
\label{A8}
\end{equation}
with eigenvalue $\lambda = 1$. Clearly, $\sigma_c(L) = [0, \infty)$, thus the eigenvalue $\lambda = 1$ is embedded in the continuous spectrum of $L$.

\medskip

\textbf{Remark 1.} A similar example can be constructed if, instead of $\theta(x)$ of \eqref{A1}, we start with the function
\begin{equation*}
\tilde{\theta}(x) =
\left\{
  \begin{array}{ccc}
    \cos x, & \quad &\ 0 \leq x < \pi/4; \\
    \\
    (e^{\pi/4} / \sqrt{2}) \, e^{-x}, & \quad &\ \pi/4 \leq x,  \\
  \end{array}
\right.
\end{equation*}
and then, in place of $f(x)$ above, we use the even extension, say $\tilde{f}(x)$, of $\tilde{\theta}(x)$ on $\mathbb{R}$.

\section{Piecewise constant potentials}
Here we show how to construct piecewise constant compactly supported potentials $q(x)$ with positive eigenvalues. The construction is based on the
following proposition.

\medskip

\textbf{Proposition 1.} Let
\begin{equation}
Q_0(x) =
\left\{
  \begin{array}{rrl}
     A, & \quad &\ a \leq x < b;
    \\
   \
    \\
     -B, & \quad &\ b \leq x,
      \\
  \end{array}
\right.
\label{C1}
\end{equation}
We consider the solution $u(x)$ of the equation
\begin{equation}
u''''(x) + Q_0(x) u(x) = 0,
\qquad
x > a,
\label{C2}
\end{equation}
satisfying the initial conditions
\begin{equation}
u(a) = \alpha_0 > 0,
\quad
u'(a) = \alpha_1 > 0,
\quad
u''(a) = \alpha_2 > 0,
\quad
u'''(a) = \alpha_3 > 0,
\label{C3}
\end{equation}
where $\alpha_0, \alpha_1, \alpha_2, \alpha_3$ are fixed. Then, given $a$, we can choose $b > a$, $A > 0$, and $B > 0$ so that
\begin{equation}
u'(\zeta) = u'''(\zeta) = 0,
\label{C4}
\end{equation}
for some $\zeta \in (b, \infty)$.

\medskip

Before giving the proof of Proposition 1, let us see how it helps to construct our desired example.

Assume without loss of generality (by shifting to the left by $\zeta$) that $\zeta = 0$ (thus, $a < b < 0$) and consider the equation
\begin{equation}
u''''(x) + q(x) u(x) = k_0^4 u(x),
\qquad
x \in \mathbb{R},
\label{C5}
\end{equation}
where $k_0 > 0$ is fixed and
\begin{equation}
q(x) =
\left\{
  \begin{array}{rrl}
     0, & \quad &\ x < a;
    \\
   Q_0(x) + k_0^4, & \quad &\ a \leq x < 0;
    \\
     q(-x), & \quad &\ 0 \leq x,
      \\
  \end{array}
\right.
\label{C6}
\end{equation}
$Q_0(x)$ being as in Proposition 1. Next, we choose the initial conditions of \eqref{C3} to be
\begin{equation}
\alpha_0 = e^{k_0 a},
\quad
\alpha_1 = k_0 e^{k_0 a},
\quad
 \alpha_2 = k_0^2 e^{k_0 a},
\quad\text{and}\quad
\alpha_3  = k_0^3 e^{k_0 a}
\label{C7}
\end{equation}
(notice that they are all positive) and consider the (unique) solution $g(x)$ of \eqref{C5} satisfying
\begin{equation}
g(a) = e^{k_0 a},
\quad
g'(a) = k_0 e^{k_0 a},
\quad
g''(a) = k_0^2 e^{k_0 a},
\quad\text{and}\quad
g'''(a)  = k_0^3 e^{k_0 a}.
\label{C7a}
\end{equation}
In other words, $g(x)$ is the unique solution of \eqref{C5} such that
\begin{equation}
g(x) = e^{k_0 x}
\qquad\text{for }\;
x < a.
\label{C8}
\end{equation}
Then, by Proposition 1 we know that
$g(x)$ satisfies \eqref{C4} (with $\zeta = 0$), i.e.
\begin{equation}
g'(0) = g'''(0) = 0,
\label{C9}
\end{equation}
while from \eqref{C6} we see that our $q(x)$ is even. Therefore, by \eqref{C9} we get that $g(x)$, too, is even and, finally, by \eqref{C8} we
can conclude that $g(x)$ is an $L_2(\mathbb{R})$-eigenfunction of the operator
\begin{equation}
Lu = u'''' + q(x) u
\label{C10}
\end{equation}
with eigenvalue $\lambda = k_0^4 > 0$.

\subsection{Proof of Proposition 1}
We start with an observation.

\medskip

\textbf{Observation 1.} Suppose that for some constant $A > 0$ we have
\begin{equation}
u''''(x) = -A u(x)
\qquad
x \in \mathbb{R},
\label{B2}
\end{equation}
with
\begin{equation}
u(a) > 0,
\qquad
u'(a) > 0,
\qquad
u''(a) > 0,
\quad\text{and}\quad
u'''(a) > 0,
\label{B3}
\end{equation}
for some fixed $a \in \mathbb{R}$. If for $j = 0, 1, 2, 3$ we denote by $x_j$ the smallest zero of $u^{(j)}(x)$ in the interval $(a, \infty)$, then
it is not hard to see (since for $u^{(j)}(x)$ to become negative, $j = 0, 1, 2$, we need first $u^{(j+1)}(x)$ to become negative) that
\begin{equation}
a < x_3 < x_2 < x_1 < x_0 < \infty.
\label{B4}
\end{equation}

\medskip

Let us, next, consider the initial value problem
\begin{equation}
u''''(x) = B u(x),
\qquad
x \in \mathbb{R},
\label{B5}
\end{equation}
\begin{equation}
u(0) = \gamma_0,
\quad
u'(0) = \gamma_1,
\quad
u''(0) = \gamma_2,
\quad
u'''(0) = \gamma_3,
\label{B6}
\end{equation}
where $\gamma_0, \gamma_1, \gamma_2, \gamma_3$ are fixed real numbers and $B$ is a (moving) parameter.

For the purpose of our analysis we need to derive a formula for the (partial) derivative of the solution $u(x) = u(x; B)$ of the above problem with
respect to the parameter $B$.

Differentiating \eqref{B5} and \eqref{B6} with respect to $B$ yields
\begin{equation}
\partial_B u''''(x; B) = B\, \partial_B u(x; B) + u(x; B),
\qquad
x \in \mathbb{R},
\label{B7}
\end{equation}
\begin{equation}
\partial_B u(0; B) = \partial_B u'(0; B) = \partial_B u''(0; B) = \partial_B u'''(0; B) = 0,
\label{B7a}
\end{equation}
where primes will always denote derivatives with respect to $x$.

It is not hard to check that the solution $\partial_B u(x; B)$ of the initial value problem \eqref{B7}--\eqref{B7a} is
\begin{equation}
\partial_B u(x; B) = \int_0^x U(x - \xi; B)\, u(\xi; B) \,d\xi,
\label{D1}
\end{equation}
where $U(x) = U(x; B)$ is the solution of the initial value problem
\begin{equation}
U''''(x) = B \, U(x),
\qquad
x \in \mathbb{R},
\label{D2}
\end{equation}
\begin{equation}
U(0) = 0,
\quad
U'(0) = 0,
\quad
U''(0) = 0,
\quad
U'''(0) = 1,
\label{D3}
\end{equation}
namely
\begin{equation}
U(x; B) = \frac{1}{2 B^{3/4}} \sinh\left(B^{1/4} x\right) - \frac{1}{2 B^{3/4}} \sin\left(B^{1/4} x\right).
\label{D4}
\end{equation}
Notice that $U(x; B)$ is entire in $B$ and positive for $x, B > 0$, and the same is true for its $x$-derivatives
\begin{equation}
U'(x; B) = \frac{1}{2 B^{1/2}} \cosh\left(B^{1/4} x\right) - \frac{1}{2 B^{1/2}} \cos\left(B^{1/4} x\right),
\label{D4a}
\end{equation}
\begin{equation}
U''(x; B) = \frac{1}{2 B^{1/4}} \sinh\left(B^{1/4} x\right) + \frac{1}{2 B^{1/4}} \sin\left(B^{1/4} x\right),
\label{D4b}
\end{equation}
and
\begin{equation}
U'''(x; B) = \frac{1}{2} \cosh\left(B^{1/4} x\right) + \frac{1}{2} \cos\left(B^{1/4} x\right).
\label{D4c}
\end{equation}

We are now ready for the main ingredient of the proof of Proposition 1.

\medskip

\textbf{Lemma 1.} Consider the equation
\begin{equation}
u''''(x) = B u(x),
\qquad
x \in \mathbb{R},
\label{B8}
\end{equation}
where $B > 0$ is a parameter and let $u(x)$ be the solution of \eqref{B8} satisfying
\begin{equation}
u(b) = \gamma_0 > 0,
\quad
u'(b) = \gamma_1 > 0,
\quad
u''(b) = \gamma_2 < 0,
\quad
u'''(b) = \gamma_3 < 0,
\label{B9}
\end{equation}
where the numbers $b, \gamma_0, \gamma_1, \gamma_2, \gamma_3$ are fixed, hence independent of $B$.
If for $j = 0, 1, 2, 3$ we denote by $z_j$ the smallest zero of $u^{(j)}(x)$ in the interval $(b, \infty)$, with the convention that $z_j = \infty$
in the case where $u^{(j)}(x)$ does not vanish in $(b, \infty)$, then we can choose $B$ so that
\begin{equation}
z_1 = z_3.
\label{B10}
\end{equation}

\smallskip

\textit{Proof}. Before we start, let us notice that without loss of generality we can take
\begin{equation*}
b = 0.
\end{equation*}

The proof is divided into three parts: In the first part we determine the behavior of $u^{(j)}(x)$, $j = 0, 1, 2, 3$,
on $x \in [0, \infty)$, in the case
where $B$ takes large values, while in the second part we determine the behavior of the same quantities in the case where $B$ is close to $0$.
Finally, in the third part of the proof we show that as $B$ moves (continuously) from larger to smaller values, it hits a value for which \eqref{B10} is achieved.

\smallskip

\textbf{I.} The assumption $u'(0) = \gamma_1 > 0$ implies that $u(x)$ is increasing on $[0, z_1)$; in particular,
\begin{equation}
u(x) \geq u(0) = \gamma_0,
\qquad
x \in [0, z_1).
\label{B11}
\end{equation}
If we integrate \eqref{B8} repeatedly (with respect to $x$) and use \eqref{B9} and \eqref{B11}, we obtain
\begin{equation}
u'''(x) \geq \gamma_3 + B \gamma_0 x,
\qquad
x \in [0, z_1),
\label{B11a}
\end{equation}
\begin{equation}
u''(x) \geq \gamma_2 + \gamma_3 x + \frac{B \gamma_0}{2} x^2,
\qquad
x \in [0, z_1),
\label{B11b}
\end{equation}
and
\begin{equation}
u'(x) \geq \gamma_1 + \gamma_2 x + \frac{\gamma_3}{2} x^2 + \frac{B \gamma_0}{6} x^3,
\qquad
x \in [0, z_1).
\label{B11c}
\end{equation}
Since $\gamma_0, \gamma_1 > 0$, it is clear that there is $\delta > 0$ and a $B^{\sharp} > 0$ such that the right-hand side of \eqref{B11c} is
$\geq \delta$ for every $B > B^{\sharp}$. In particular, if $z_1 < \infty$, then we would have from \eqref{B11c} that $u'(z_1) \geq \delta$, which is impossible since $u'(z_1) = 0$. Therefore $z_1 = \infty$ for every $B > B^{\sharp}$, i.e.
\begin{equation}
\text{if }\; B > B^{\sharp},
\qquad\text{then}\qquad
 u'(x) > 0
\quad\text{for all }\;
x \in [0, \infty).
\label{B12}
\end{equation}
From \eqref{B12} we get that $u(x)$ is increasing and, hence $u(x) \geq \gamma_0 > 0$ on $[0, \infty)$. Thus \eqref{B8} implies that $u'''(x)$ is
(strictly) increasing and eventually positive on $[0, \infty)$. Consequently, $u'''(x)$ has a unique zero $z_3 \in (0, \infty)$. Finally, $u''(x)$
has a unique minimum on $[0, \infty)$ attained at $x = z_3$; as for $x > z_3$, the function $u''(x)$ is (strictly) increasing and eventually
positive, thus $z_3 < z_2 < \infty$.

\smallskip

\textbf{II.}  The solution $u(x) = u(x; B)$ of \eqref{B8} depends analytically in the parameter $B$; actually, it is entire in $B$.
Thus, for every $x \in \mathbb{R}$ and every $B \in \mathbb{C}$ we have
\begin{equation}
u(x; B) = \sum_{n=0}^{\infty} u_n(x) B^n,
\label{B13}
\end{equation}
where
\begin{equation}
u_n(x) = \left.\frac{\partial^n u(x; B)}{\partial B^n} \right|_{B=0},
\label{B13a}
\end{equation}
in particular
\begin{equation}
u_0(x) = u(x; 0) = \gamma_0 + \gamma_1 x + \frac{\gamma_2}{2} x^2 + \frac{\gamma_3}{6} x^3.
\label{B13b}
\end{equation}
Differentiating \eqref{B8} $n$ times with respect to $B$ and then setting $B=0$ yields
\begin{equation}
u_n''''(x) = n u_{n-1}(x)
\qquad
x \in \mathbb{R},
\quad
n = 1, 2, \ldots
\label{B13c}
\end{equation}
(as usual, primes denote derivatives with respect to $x$), with initial conditions (due to \eqref{B9})
\begin{equation}
u_n(0) = u_n'(0) = u_n''(0) = u_n'''(0) = 0,
\qquad
n = 1, 2, \ldots \,.
\label{B13d}
\end{equation}
By using \eqref{B8} in \eqref{B13} we obtain
\begin{equation}
u''''(x; B) = B u(x; B) = \sum_{n=0}^{\infty} u_n(x) B^{n+1},
\label{B13e}
\end{equation}
thus, by integrating \eqref{B13e} with respect to $x$ repeatedly, we can obtain expansions in powers of $B$, similar to \eqref{B13} for the $x$-derivatives $u'(x; B)$, $u''(x; B)$, and $u'''(x; B)$ (which are, therefore, entire in $B$ too).

From the analytic dependence and the form of $u(x; 0)$, as given in \eqref{B13b}, it is not hard to see that there is a $B^{\flat} > 0$ such that
if $0 < B < B^{\flat}$, then: (i) $u'''(x; B)$ is negative for $x \in [0, \infty)$, with $\lim_{x \to \infty} u'''(x; B) = -\infty$;
(ii) $u''(x; B)$ is negative and decreasing in $x$ on $[0, \infty)$, with $\lim_{x \to \infty} u''(x; B) = -\infty$;
(iii) $u'(x; B)$ is decreasing in $x$ on $[0, \infty)$, with $\lim_{x \to \infty} u'(x; B) = -\infty$, thus $z_1 < \infty$;
and (iv) $u(x; B)$ has a unique maximum in $[0, \infty)$ (attained at $x = z_1$), with $\lim_{x \to \infty} u(x; B) = -\infty$, hence
$z_1 <z_0 < \infty$.

\smallskip

\textbf{III.} Suppose we start with a value of $B$ larger than $B^{\sharp}$ and we decrease it towards $0$ in a continuous motion. As we have seen in Parts I and II, $u'(x)$ has a unique minimum in $[0, \infty)$ attained at $x = z_2$. Initially $u'(z_2) > 0$, but as $B$ decreases it will reach a
value, say $B_1$ for which $z_1 = z_2$, so that $u'(z_1) = u''(z_1) = u''(z_2) = 0$. In addition, we must have that $u'''(z_2) > 0$ since $u''(x)$
is negative and decreasing near $x = 0$ and then starts increasing and becomes $0$ at $x = z_2$; thus $u''(x)$ attains its minimum at
$x = z_3 < z_2$, which implies that $u'''(z_2) > 0$.

Therefore, if $B$ becomes slightly smaller than $B_1$, then we will have $z_3 < z_1 < z_2 < \infty$. Let us now examine the quantities $dz_1/dB$
and $dz_3/dB$. Since $u'(z_1; B) = 0$ and $u'''(z_3; B) = 0$, by implicit differentiation we get
\begin{equation}
u''(z_1; B) \frac{dz_1}{dB} = - \partial_B u'(z_1; B)
\label{B14}
\end{equation}
and
\begin{equation}
u''''(z_3; B) \frac{dz_1}{dB} = - \partial_B u'''(z_3; B)
\quad\text{or}\quad
B u(z_3; B) \frac{dz_1}{dB} = - \partial_B u'''(z_3; B)
\label{B15}
\end{equation}
(since $u'''' = Bu$).
Now
\begin{equation}
u''(z_1; B) < 0
\qquad\text{and}\qquad
u(z_3; B) > 0.
\label{B16}
\end{equation}
Also, by differentiating formula \eqref{D1} and using \eqref{D3} we get
\begin{equation}
\partial_B u'(z_1; B) = \int_0^{z_1} U'(z_1 - \xi; B)\, u(\xi; B) \,d\xi
\label{B17}
\end{equation}
and
\begin{equation}
\partial_B u'''(z_3; B) = \int_0^{z_3} U'''(z_3 - \xi; B)\, u(\xi; B) \,d\xi
\label{B18}
\end{equation}
In view \eqref{D4a} and \eqref{D4c} the quantity $U'(z_1 - \xi; B)$ is (strictly) positive for $\xi \in (0, z_1)$ and $U'''(z_3 - \xi; B)$
is (strictly) positive for $\xi \in (0, z_3)$. Also, $u(\xi; B) > 0$ for $\xi \in (0, z_1)$. Hence, \eqref{B17} and \eqref{B18} imply that
\begin{equation}
\partial_B u'(z_1; B) > 0
\qquad\text{and}\qquad
\partial_B u'''(z_3; B) > 0.
\label{B19}
\end{equation}
Using \eqref{B16} and \eqref{B19} in \eqref{B14} and \eqref{B15} we can conclude that
\begin{equation}
\frac{dz_1}{dB} > 0,
\qquad\text{while}\qquad
\frac{dz_3}{dB} < 0,
\label{B20}
\end{equation}
which tells us that $z_1$ decreases, while $z_3$ increases, as $B$ decreases.
As we have seen, for $B$ slightly smaller than $B_1$ we have that $z_3 < z_1$. Also, as $B$ decreases $z_3$, remains smaller than $z_0$ as long as it exists as a real number, and ceases to exist at the moment when $z_3 = z_0$ (since $z_0$ is the smallest positive zero of $u'''' = B u$). At the moment, though, when $z_3 = z_0$ we must have $z_1 < z_0$. Therefore there must be a $B > 0$ for which $z_1(B) = z_3(B)$.
\hfill $\blacksquare$

\medskip

We are now ready to finish the proof of Proposition 1.

\medskip

\textit{Proof of Proposition 1}. Given $a$ we fix an $A > 0$. Then we choose a $b \in (x_2, x_1)$, where $x_2, x_1$ are as in \eqref{B4} of Observation 1. Clearly, $b > a$ and for this $b$ it is easy to see that \eqref{B9} is satisfied. Finally, from Lemma 1 we obtain the existence of a
$B > 0$ for which $z_1 = z_3$, which means that there is a positive $B$ for which \eqref{C4} holds.
\hfill $\blacksquare$

\section{The square of a Schr\"{o}dinger operator}
Finally, let us mention a somehow trivial example where the fourth-order operator $L$ is not of the form $d^4/dx^4 + q(x)$.

Let $H = -d^2/dx^2 + V(x)$ be a Schr\"{o}dinger operator, acting on $L_2(\mathbb{R})$, whose potential $V(x)$ is smooth, say $C^r$, $r \geq 2$.
Suppose that $V(x) \to 0$ as $x \to \pm\infty$ and that $H$ has bound states $\kappa_1 < \kappa_2 < \cdots < 0$. Then, the operator
\begin{equation}
L = H^2 = \frac{d^4}{dx^4} - V(x) \frac{d^2}{dx^2} - 2 V'(x) \frac{d}{dx} + V(x)^2 - V''(x)
\label{S2}
\end{equation}
has continuous spectrum $\sigma_c(L) = [0, \infty)$ and eigenvalues $\kappa_1^2 > \kappa_2^2 >  \cdots > 0$.

\end{document}